\documentclass[11pt,twoside,a4paper]{article}
\usepackage{amsmath,amssymb,amsthm}

\usepackage{graphicx,psfrag}
\newtheorem{thm}{Theorem}[section]
\newtheorem{lem}[thm]{Lemma}
\newtheorem{pro}[thm]{Proposition}
\newtheorem{cor}[thm]{Corollary}

\theoremstyle{definition}

\theoremstyle{remark}
\newtheorem{rem}[thm]{Remark}

\newcommand{\R}{\mathbb{R}}
\newcommand{\Z}{\mathbb{Z}}
\newcommand{\N}{\mathbb{N}}
\newcommand{\C}{\mathbb{C}}
\newcommand{\Q}{\mathbb{Q}}

\newcommand{\cA}{\mathcal{A}}

\newcommand{\cE}{\mathcal{E}}

\newcommand{\cP}{\mathcal{P}}

\newcommand{\cU}{\mathcal{U}}
\newcommand{\cV}{\mathcal{V}}
\newcommand{\cW}{\mathcal{W}}

\newcommand{\de}{\delta}

\newcommand{\ep}{\varepsilon}

\newcommand{\si}{\sigma}

\newcommand{\la}{\lambda}
\newcommand{\La}{\Lambda}
\renewcommand{\phi}{\varphi}

\newcommand{\dist}{\operatorname{dist}}
\newcommand{\diam}{\operatorname{diam}}

\newcommand{\hyp}{\operatorname{H}}

\newcommand{\asdim}{\operatorname{asdim}}
\newcommand{\hypdim}{\operatorname{hypdim}}

\newcommand{\cone}{\operatorname{Co}}
\newcommand{\pt}{\operatorname{pt}}
\newcommand{\mesh}{\operatorname{mesh}}
\newcommand{\cdim}{\operatorname{cdim}}
\newcommand{\Cdim}{\operatorname{Cdim}}

\newcommand{\an}{\operatorname{An}}

\newcommand{\es}{\emptyset}
\renewcommand{\d}{\partial}
\newcommand{\di}{\d_{\infty}}
\newcommand{\set}[2]{\{#1:\,\text{#2}\}}
\newcommand{\sm}{\setminus}
\newcommand{\sub}{\subset}

\newcommand{\ov}{\overline}
\newcommand{\wt}{\widetilde}
\newcommand{\wh}{\widehat}

\begin{document}

\title{Dimensions of locally and asymptotically
self-similar spaces}
\author{Sergei Buyalo
and Nina Lebedeva\footnote{Both authors are supported by RFFI Grant
05-01-00939 and Grant NSH-1914.2003.1}}

\date{}
\maketitle

\begin{abstract} We obtain two in a sense dual to each other
results: First, that the capacity dimension of every compact,
locally self-similar metric space coincides
with the topological dimension, and second, that the asymptotic
dimension of a metric space, which is asymptotically similar
to its compact subspace coincides
with the topological dimension of the subspace.
As an application of the first result, we prove the Gromov conjecture
that the asymptotic dimension of every hyperbolic group
$G$
equals the topological
dimension of its boundary at infinity plus 1,
$\asdim G=\dim\di G+1$.
As an application of the second result, we construct Pontryagin
surfaces for the asymptotic dimension, in particular, those are
first examples of metric spaces
$X$, $Y$
with
$\asdim(X\times Y)<\asdim X+\asdim Y$.
Other applications are also given.
\end{abstract}

\tableofcontents

\section{Introduction}

We say that a map
$f:Z\to Z'$
between metric spaces is {\em quasi-homothetic}
with coefficient
$R>0$,
if for some
$\la\ge 1$
and for all
$z$, $z'\in Z$,
we have
$$R|zz'|/\la\le|f(z)f(z')|\le\la R|zz'|.$$
In this case, we also say that
$f$
is
$\la$-quasi-homothetic
with coefficient
$R$.

A metric space
$Z$
is {\em locally similar} to a metric space
$Y$,
if there is
$\la\ge 1$
such that for every sufficiently large
$R>1$
and every
$A\sub Z$
with
$\diam A\le\La_0/R$,
where
$\La_0=\min\{1,\diam Y/\la\}$,
there is a
$\la$-quasi-homothetic
map
$f:A\to Y$
with coefficient
$R$
(note that the condition
$\diam A\le\La_0/R$
implies
$\diam f(A)\le\diam Y$).
If a metric space
$Z$
is locally similar to itself then we say that
$Z$
is {\em locally self-similar}.

The notion of the capacity dimension of a metric space
$Z$, $\cdim Z$,
is introduced in \cite{Bu1}, and turns out to be useful
in many questions, \cite{Bu2}. The capacity dimension
is larger than or equal to the topological dimension,
$\dim Z\le\cdim Z$
for every metric space
$Z$,
and it is important to know for which spaces the equality
holds. Our first main result is the
following

\begin{thm}\label{thm:main} Assume that a metric space
$Z$
is locally similar to a compact metric space
$Y$.
Then
$\cdim Z<\infty$
and
$\cdim Z\le\dim Y$.
\end{thm}

\begin{cor}\label{cor:selfsim} The capacity dimension
of every compact, locally self-simi\-lar metric space
$Z$
is finite and coincides with its topological dimension,
$\cdim Z=\dim Z$.
\end{cor}

In contrast, we also prove a proposition (see
Proposition~\ref{pro:locsimcdim}), which allows to
construct examples of compact metric spaces with the
capacity dimension arbitrarily larger than the
topological dimension.

Now, consider in a sense the dual situation. A metric space
$X$
is {\em asymptotically similar} to a metric space
$Y$,
if there are
$\La_0$, $\la\ge 1$
such that for every sufficiently large
$R>1$
and every
$A\sub X$
with
$\diam A\le R/\La_0$
there is a
$\la$-quasi-homothetic
map
$f:Y\to X$
with
coefficient
$R$,
whose image contains an isometric copy
$A'\sub X$
of
$A$, $A'\sub f(Y)$.
If a metric space
$X$
is asymptotically similar to a bounded subset then we say that
$X$
is {\em asymptotically self-similar}. Taking a copy
$A'$
instead of
$A$
provides an additional flexibility of this definition,
which is necessary for applications, see
sect.~\ref{subsect:asympont}.

We recall the well established notion of the asymptotic
dimension,
$\asdim$,
in sect.~\ref{subsect:defcapdim}.

Our second main result is the following

\begin{thm}\label{thm:main2} Assume that a metric space
$X$
is asymptotically similar to a compact metric space
$Y$.
Then the both dimensions,
$\asdim X$, $\dim Y$
are finite and coincide,
$\asdim X=\dim Y$.
\end{thm}

As applications of Theorem~\ref{thm:main}, we obtain

\begin{thm}\label{thm:cdimhypgroup} The capacity dimension of the boundary
at infinity of any hyperbolic group
$G$
(taken with any visual metric) coincides with the topological
dimension,
$\cdim\di G=\dim\di G$.
\end{thm}

Theorem~\ref{thm:cdimhypgroup} together with the main result of
\cite{Bu1} leads to the following result which proves a
Gromov conjecture, see \cite[1.$\text{E}_1'$]{Gr}.

\begin{thm}\label{thm:gromovconjecture} The asymptotic dimension
of any hyperbolic group
$G$
equals topological dimension
of its boundary at infinity plus 1,
$$\asdim G=\dim\di G+1.$$
\end{thm}

Another application of Theorem~\ref{thm:cdimhypgroup} is
the following embedding result, obtained from the main
result of \cite{Bu2}.

\begin{thm}\label{thm:embd} Every hyperbolic group
$G$
admits a quasi-isometric embedding
$G\to T_1\times\dots\times T_n$
into the
$n$-fold
product of simplicial metric trees
$T_1,\dots,T_n$
with
$n=\dim\di G+1$.
\end{thm}

The group structure plays no role in the proof of
Theorems~\ref{thm:cdimhypgroup} -- \ref{thm:embd}.
Actually, we prove more general Theorems~\ref{thm:cdimhypspace},
\ref{thm:asdim}, \ref{thm:embdspace}, and have chosen
the statements above for simplicity of formulations.

Theorem~\ref{thm:cdimhypspace} has applications also to
nonembedding results, which are discussed in
sect.~\ref{subsect:embnonemb}, see Theorem~\ref{thm:nonembd}.
Using Corollary~\ref{cor:selfsim}, we give
examples of strict inequality in the product theorem for
the capacity dimension. These are famous Pontryagin
surfaces self-similar construction of which is discussed
in sect.~\ref{sect:sspont}, see Theorem~\ref{thm:sspont}.
Finally, we construct metric spaces asymptotically
similar to self-similar Pontryagin surfaces.
As a corollary of Theorem~\ref{thm:main2}, we give
examples of strict inequality in the product theorem
for the asymptotic dimension, that is, we construct
metric spaces
$X$, $Y$
with
$$\asdim(X\times Y)<\asdim X+\asdim Y,$$
see Corollary~\ref{cor:asympont}.

{\em Acknowledgements.} The first author is grateful to
Bruce Kleiner and Viktor Schroeder for numerous
stimulating discussions related to this paper.
We are also grateful to Alex Dranishnikov for
informing us about the paper \cite{Sw}.

\section{Preliminaries}

Here, we recall notions and facts necessary for the paper.

Let
$Z$
be a metric space. For
$U$, $U'\sub Z$
we denote by
$\dist(U,U')$
the distance between
$U$
and
$U'$,
$\dist(U,U')=\inf\set{|uu'|}{$u\in U,\ u'\in U'$}$
where
$|uu'|$
is the distance between
$u$, $u'$.
For
$r>0$
we denote by
$B_r(U)$
the open
$r$-neighborhood
of
$U$, $B_r(U)=\set{z\in Z}{$\dist(z,U)<r$}$,
and by
$\ov B_r(U)$
the closed
$r$-neighborhood
of
$U$, $\ov B_r(U)=\set{z\in Z}{$\dist(z,U)\le r$}$.
We extend these notations over all real
$r$
putting
$B_r(U)=U$
for
$r=0$,
and defining
$B_r(U)$
for
$r<0$
as the complement of the closed
$|r|$-neighborhood
of
$Z\sm U$,
$B_r(U)=Z\sm\ov B_{|r|}(Z\sm U)$.

Given a family
$\cU$
of subsets in a metric space
$Z$
we define
$\mesh(\cU)=\sup\set{\diam U}{$U\in\cU$}$.
The {\em multiplicity} of
$\cU$, $m(\cU)$,
is the maximal number of members of
$\cU$
with nonempty intersection. We say that a family
$\cU$
is {\em disjoint} if
$m(\cU)=1$.

A family
$\cU$
is called a {\em covering} of
$Z$
if
$\cup\set{U}{$U\in\cU$}=Z$.
A covering
$\cU$
is said to be {\em colored} if it is the union
of
$m\ge 1$
disjoint families,
$\cU=\cup_{a\in A}\cU^a$, $|A|=m$.
In this case we also say that
$\cU$
is
$m$-colored.
Clearly, the multiplicity of a
$m$-colored
covering is at most
$m$.

Let
$\cU$
be a family of open subsets in a metric space
$Z$
which cover
$A\sub Z$.
Given
$z\in A$,
we let
$$L(\cU,z)=\sup\set{\dist(z,Z\sm U)}{$U\in\cU$}$$
be the Lebesgue number of
$\cU$
at
$z$, $L(\cU)=\inf_{z\in A}L(\cU,z)$
be the Lebesgue number
of the covering
$\cU$
of
$A$.
For every
$z\in A$,
the open ball
$B_r(z)\sub Z$
of radius
$r=L(\cU)$
centered at
$z$
is contained in some member of the covering
$\cU$.

We shall use the following obvious fact (see e.g. \cite{Bu1}).

\begin{lem}\label{lem:insidecov} Let
$\cU$
be an open covering of
$A\sub Z$
with
$L(\cU)>0$.
Then for every
$s\in(0,L(\cU))$
the family
$\cU_{-s}=B_{-s}(\cU)$
is still an open covering of
$A$.
\qed
\end{lem}

\subsection{Definitions of capacity and
asymptotic dimensions}\label{subsect:defcapdim}

There are several equivalent definitions of the capacity
dimension, see \cite{Bu1}. In this paper, we use the
following two.

(i) The capacity dimension of a metric space
$Z$, $\cdim(Z)$,
is the minimal integer
$m\ge 0$
with the following property: There is a constant
$\de\in(0,1)$
such that for every sufficiently small
$\tau>0$
there exists a
$(m+1)$-colored
open covering
$\cU$
of
$Z$
with
$\mesh(\cU)\le\tau$
and
$L(\cU)\ge\de\tau$.

(ii) The capacity dimension of a metric space
$Z$, $\cdim(Z)$,
is the minimal integer
$m\ge 0$
with the following property: There is a constant
$\de\in(0,1)$
such that for every sufficiently small
$\tau>0$
there exists an open covering
$\cU$
of
$Z$
with multiplicity
$m(\cU)\le m+1$,
for which
$\mesh(\cU)\le\tau$
and
$L(\cU)\ge\de\tau$.

The asymptotic dimension is a quasi-isometry invariant of
a metric space introduced in \cite{Gr}.
There are also several equivalent definitions,
see \cite{Gr}, \cite{BD}, and we use the following one.
The asymptotic dimension of a metric space
$X$
is the minimal integer
$\asdim X=n$
such that for every positive
$d$
there is an open covering
$\cU$
of
$X$
with
$m(\cU)\le n+1$, $\mesh(\cU)<\infty$
and
$L(\cU)\ge d$.

The following notion turns out to be useful for our
purposes. This notion is called the Higson property in
\cite[sect.~4]{DZ} and the asymptotic dimension of linear
type in the book \cite[Example~9.14]{Ro}.
We call it {\em asymptotic capacity dimension}. Similarly
to the capacity dimension, the following two definitions,
colored and covering ones, are equivalent.

(i) The asymptotic capacity dimension of a metric space
$X$, $\Cdim(X)$,
is the minimal integer
$m\ge 0$
with the following property: There is a constant
$\de\in(0,1)$
such that for every large
$R>1$
there exists a
$(m+1)$-colored
open covering
$\cU$
of
$X$
with
$\mesh(\cU)\le R$
and
$L(\cU)\ge\de R$.

(ii) The asymptotic capacity dimension of a metric space
$X$, $\Cdim(X)$,
is the minimal integer
$m\ge 0$
with the following property: There is a constant
$\de\in(0,1)$
such that for every sufficiently large
$R>1$
there exists an open covering
$\cU$
of
$X$
with multiplicity
$m(\cU)\le m+1$,
for which
$\mesh(\cU)\le R$
and
$L(\cU)\ge\de R$.

Clearly,
$\asdim X\le\Cdim X$
for any metric space
$X$.

\section{Auxiliary facts}

Here, we collect some facts needed for the proof of
our main results.

The following lemma implies in particular the finite
union theorem for the capacity dimension, which is similar
to the appropriate theorems for the asymptotic dimension
\cite{BD} and the Assouad-Nagata dimension \cite{LS}.

\begin{lem}\label{lem:union} Suppose, that
$Z$
is a metric space and
$A$, $B\sub Z$.
Let
$\cU$
be an open covering of
$A$,
$\cV$
be an open covering of
$B$
both with multiplicity at most
$m$.
If
$\mesh(\cV)\le L(\cU)/2$
then there exist an open covering
$\cW$
of
$A\cup B$
with multiplicity at most
$m$
and
$\mesh(\cW)\le\max\{\mesh(\cV),\mesh(\cU)\}$,
$L(\cW)\ge\min\{L(\cU)/2,L(\cV)\}$.
\end{lem}

\begin{proof} We can assume that
$L(\cU)<\infty$,
i.e. no member of
$\cU$
covers
$Z$,
because otherwise we take
$\cW=\cU$.

We put
$r=L(\cU)/2$
and consider the family
$\wt\cU=\set{B_{-r}(U)}{$U\in\cU$}$.
This family still covers
$A$.
Next, let
$\wt\cV$
be the family of all
$V\in\cV$,
each of which intersects some
$\wt U\in\wt\cU$.
For every
$\wt V\in\wt\cV$,
we fix
$\wt U\in\wt\cU$
with
$\wt U\cap\wt V\neq\es$
and consider the union
$W=W(\wt U)$
of
$\wt U$
and all
$\wt V\in\wt\cV$
assigned in this way to
$\wt U$.
Now, we take the family
$\cW$
consisting of all
$V\in\cV$,
which do not enter
$\wt\cV$,
and all
$W(\wt U)$, $\wt U\in\wt\cU$.

By the remark at the beginning, the family
$\cW$
covers
$A\cup B$.
Since
$\diam V\le r$
for every
$V\in\cV$,
we have
$W(\wt U)\sub U$
for the corresponding
$U\in\cU$
and thus
$\diam W(\wt U)\le\mesh\cU$.
Hence,
$\mesh\cW\le\max\{\mesh(\cU),\mesh\cV\}$.

Clearly, for the Lebesgue number of
$\cW$
we have
$$L(\cW)\ge\min\{L(\cV),L(\wt\cU)\}\ge
   \min\{L(\cV),L(\cU)/2\}.$$
Finally, let
$\cA$
be a collection of members of
$\cW$
with nonempty intersection. If
$\cA$
contains no member of
$\cV\sm\wt\cV$
then every member of
$\cA$
is contained in some member of
$\cU$
and thus
$|\cA|\le m$.
Otherwise, the intersection of
$\cA$
is contained in some
$V\in\cV\sm\wt\cV$
and therefore, it misses the closure of any
$\wt U\in\wt\cU$.
Thus,
$\cA$
contains at most as many members as the number of
members of
$\cV$
contains the intersection. This shows that the multiplicity of
$\cW$
is at most
$m$.
\end{proof}

We say that a metric space
$Z$
is {\em doubling at small scales} if there is a constant
$N\in\N$
such that for every sufficiently small
$r>0$
every ball in
$Z$
of radius
$2r$
can be covered by at most
$N$
balls of radius
$r$.

Similarly, a metric space
$X$
is {\em asymptotically doubling} if there is a constant
$N\in\N$
such that for every sufficiently large
$R>1$
every ball in
$X$
of radius
$2R$
can be covered by at most
$N$
balls of radius
$R$.

\begin{lem}\label{lem:scaledoubling} (1) Assume that a metric
space
$Z$
is locally similar to a compact metric space
$Y$.
Then
$Z$
is doubling at small scales.

(2) Assume that a metric space
$X$
is asymptotically similar to a compact metric space
$Y$.
Then
$X$
is asymptotically doubling.
\end{lem}

\begin{proof} (1) There is
$\la\ge 1$
such that for every sufficiently large
$R>1$
and every
$A\sub Z$
with
$\diam A\le\La_0/R$, $\La_0=\min\{1,\diam Z/\la\}$,
there is a
$\la$-quasi-homothetic
map
$f:A\to Z$
with coefficient
$R$.

We fix a positive
$\rho\le\La_0/(4\la)$.
Because
$Y$
is compact, there is
$N\in\N$
such that any subset
$B\sub Y$
can be covered by at most
$N$
balls of radius
$\rho$
centered at points of
$B$.
Take
$r>0$
small enough so that
$R=\la\rho/r$
satisfies the assumption above. Then for any ball
$B_{2r}\sub Z$,
we have
$$\diam B_{2r}\le 4r\le\La_0/R,$$
and thus there is a
$\la$-quasi-homothetic
map
$f:B_{2r}\to Y$
with coefficient
$R$.
The image
$f(B_{2r})$
is covered by at most
$N$
balls of radius
$\rho$
centered at points of
$f(B_{2r})$.
The preimage under
$f$
of every such a ball is contained in a ball of radius
$\le\la\rho/R=r$
centered at a point in
$B_{2r}$.
Hence,
$B_{2r}$
is covered by at most
$N$
balls of radius
$r$,
and
$Z$
is doubling at small scales.

(2) Again, there are
$\La_0$, $\la\ge 1$
such that for every sufficiently large
$R>1$
and every
$A\sub X$
with
$\diam A\le R/\La_0$
there is a
$\la$-quasi-homothetic
map
$f:A\to Z$
with coefficient
$R$,
whose image contains an isometric copy
$A'\sub X$
of
$A$, $A'\sub f(Y)$.

We fix a positive
$\rho\le 1/(4\la\La_0)$.
Because
$Y$
is compact, there is
$N\in\N$
such that any subset
$B\sub Y$
can be covered by at most
$N$
balls of radius
$\rho$
centered at points of
$B$.
Take
$R>1$
large enough satisfying the assumption above. Then for
any ball
$B_{2R}\sub X$
of radius
$2R$,
there is a
$\la$-quasi-homothetic
map
$f:Y\to X$
with coefficient
$4\La_0R$,
such that
$f(Y)$
contains an isometric copy
$B_{2R}'$
of
$B_{2R}$.
Without loss of generality, we can assume that
$B_{2R}'=B_{2R}$.

Then, the preimage
$B=f^{-1}(B_{2R})$
is covered by at most
$N$
balls of radius
$\rho$
centered at points of
$B$.
The image under
$f$
of every such a ball is contained in a ball of radius
$\le 4\la\La_0\rho R\le R$
centered at a point in
$B_{2R}$.
Hence,
$B_{2R}$
is covered by at most
$N$
balls of radius
$R$,
and
$X$
is asymptotically doubling.
\end{proof}

The idea of the following lemma is borrowed from
\cite[Lemma~2.3]{LS} as well as its proof, which
we give for convenience of the reader.

\begin{lem}\label{lem:finitedim} Assume that a metric space
$Z$
is doubling at small scales, and a metric space
$X$
is asymptotically doubling. Then
$\cdim Z<\infty$
and
$\Cdim X<\infty$.
\end{lem}

\begin{proof} By the assumption, there is
$n\in\N$
such that every ball
$B_{4r}\sub Z$
of radius
$4r$
is covered by at most
$n+1$
balls
$B_{r/2}$
for all sufficiently small
$r>0$.
We fix a maximal
$r$-separated
set
$Z'\sub Z$,
i.e.
$|zz'|>r$
for each distinct
$z$, $z'\in Z'$.
Then, the family
$\cU'=\set{B_r(z)}{$z\in Z'$}$
is an open covering of
$Z$.

Since every ball
$B_{r/2}$
contains at most one point from
$Z'$
and
$B_{4r}(z)$
is covered by at most
$n+1$
balls
$B_{r/2}$,
the ball
$B_{4r}(z)$
contains at most
$n+1$
points from
$Z'$
for every
$z\in Z'$.
Thus, there is a coloring
$\chi:Z'\to A$, $|A|=n+1$,
such that
$\chi(z)\neq\chi(z')$
for each
$z$, $z'\in Z$
with
$|zz'|<4r$.

For
$a\in A$,
we let
$Z_a'=\chi^{-1}(a)$
be the set of the color
$a$.
Then
$|zz'|\ge 4r$
for distinct
$z$, $z'\in Z_a'$.
Putting
$\cU_a=\set{B_{2r}(z)}{$z\in Z_a'$}$,
we obtain an open
$(n+1)$-colored
covering
$\cU=\cup_{a\in A}\cU_a$
of
$Z$
with
$\mesh(\cU)\le 4r$
and
$L(\cU)\ge r$.
This shows that
$\cdim Z\le n$.

A similar argument shows that
$\Cdim X<\infty$.
We leave details to the reader as an exercise.
\end{proof}

We shall use the following facts obviously implied
by the definition of a quasi-homothetic map.

\begin{lem}\label{lem:distort} Let
$h:Z\to Z'$
be a
$\la$-quasi-homothetic
map with coefficient
$R$.
Let
$V\sub Z$, $\wt\cU$
be an open covering of
$h(V)$
and
$\cU=h^{-1}(\wt\cU)$.
Then
\begin{itemize}
\item[(1)] $R\mesh(\cU)/\la\le\mesh(\wt\cU)
           \le\la R\mesh(\cU)$;

\item[(2)] $\la R\cdot L(\cU)\ge L(\wt\cU)
    \ge R\cdot L(\cU)/\la$,
where
$L(\cU)$
is the Lebesgue number of
$\cU$
as a covering of
$V$.\qed
\end{itemize}
\end{lem}

\section{Proof of Theorem~\ref{thm:main}}

It follows from Lemmas~\ref{lem:scaledoubling}
and \ref{lem:finitedim}, that
$\cdim Z=N$
is finite. We can also assume that
$\dim Y=n$
is finite. There is a constant
$\de\in(0,1)$
such that for every sufficiently small
$\tau>0$
there exists a
$(N+1)$-colored
open covering
$\cV=\cup_{a\in A}\cV^a$
of
$Z$
with
$\mesh(\cV)\le\tau$
and
$L(\cV)\ge\de\tau$.
It is convenient to take
$A=\{0,\dots,N\}$
as the color set.

There is a constant
$\la\ge 1$,
such that for every sufficiently large
$R>1$
and every
$V\sub Z$
with
$\diam V\le\La_0/R$, $\La_0=\min\{1,\diam Y/\la\}$,
there is a
$\la$-quasi-homothetic
map
$h_V:V\to Y$
with coefficient
$R$.

Using that
$Y$
is compact and
$\dim Y=n$,
we find for every
$a\in A$
a finite open covering
$\wt\cU_a$
of
$Y$
with multiplicity
$m(\wt\cU_a)\le n+1$
such the following holds:

\begin{itemize}
\item[(i)] $\mesh(\wt\cU_0)\le\frac{\de}{2\la}$;
\item[(ii)] $\mesh(\wt\cU_{a+1})
    \le\frac{1}{2\la^2}\min\{L(\wt\cU_a),\mesh(\wt\cU_a)\}$
for every
$a\in A$, $a\le N-1$.
\end{itemize}
Then
$l:=\min\{L(\wt\cU_N),\frac{1}{2}L(\wt\cU_{N-1}),
   \dots,\frac{1}{2^N}L(\wt\cU_0)\}>0$
and
$\mesh(\wt\cU_a)\le\frac{\de}{2\la}$
for every
$a\in A$.

For every
$V\in\cV$,
consider the slightly smaller subset
$V'=B_{-\de\tau/2}(V)$.
Then, the sets
$Z_a=\cup_{V\in\cV^a}V'\sub Z$, $a\in A$,
cover
$Z$, $Z=\cup_{a\in A}Z_a$,
because
$L(\cV)\ge\de\tau$.

Given
$V\in\cV$,
we fix a
$\la$-quasi-homothetic
map
$h_V:V\to Z$
with coefficient
$R=1/\tau$
and put
$\wt V=h_V(V')$.
Now, for every
$a\in A$, $V\in\cV^a$
consider the family
$\wt\cU_{a,V}=\set{\wt U\in\wt\cU_a}{$\wt V\cap\wt U\neq\es$}$,
which is obviously a covering of
$\wt V$
with multiplicity
$\le n+1$.
Then,
$$\cU_{a,V}=\set{h_V^{-1}(\wt U)}{$\wt U\in\wt\cU_{a,V}$}$$
is an open covering of
$V'$
with multiplicity
$\le n+1$.

Note that
$U=h_V^{-1}(\wt U)$
is contained in
$V$
for every
$\wt U\in\wt\cU_{a,V}$
because
$\dist(v',Z\sm V)>\de\tau/2$
for every
$v'\in V'$
and
$\diam U\le\la\tau\diam\wt U\le\de\tau/2$.
Thus the family
$\cU_{a,V}$
is contained in
$V$.
Now, the family
$\cU_a=\cup_{V\in\cV^a}\cU_{a,V}$
covers the set
$Z_a$
of the color
$a$,
and it has the following properties

\begin{itemize}
\item[(1)] for every
$a\in A$,
the multiplicity of
$\cU_a$
is at most
$n+1$;

\item[(2)] $\mesh(\cU_{a+1})
   \le\frac{1}{2}\min\{L(\cU_a),\mesh(\cU_a)\}$
for every
$a\in A$, $a\le N-1$
($L(\cU_a)$
means the Lebesgue number of
$\cU_a$
as a covering of
$Z_a$);

\item[(3)] $\mesh(\cU_a)\le\la\tau\mesh(\wt\cU_a)$
and
$L(\cU_a)\ge\tau L(\wt\cU_a)/\la$
for every
$a\in A$.
\end{itemize}

Indeed, distinct
$V_1$, $V_2\in\cV^a$
are disjoint and thus any
$U_1\in\cU_{a,V_1}$, $U_2\in\cU_{a,V_2}$
are disjoint because
$U_1\sub V_1$, $U_2\sub V_2$.
This proves (1). Furthermore, for every
$a\in A$, $a\le N-1$,
and every
$U\in\cU_{a+1}$,
we have
\begin{eqnarray*}
  \diam U&\le&\la\tau\mesh(\wt\cU_{a+1})
  \le\frac{\tau}{2\la}\min\{L(\wt\cU_a),\mesh(\wt\cU_a)\}\\
  &\le&\frac{1}{2}\min\{L(\cU_a),\mesh(\cU_a)\}
\end{eqnarray*}
by Lemma~\ref{lem:distort}, hence, (2). Finally,
(3) also follows from Lemma~\ref{lem:distort}.

Now, we put
$\wh\cU_{-1}=\{Z\}$, $\wh\cU_0=\cU_0$
and assume that for some
$a\in A$,
we have already constructed families
$\wh\cU_0,\dots,\wh\cU_a$
so that
$\wh\cU_a$
covers
$Z_0\cup\dots\cup Z_a$
with multiplicity
$\le n+1$
and
$\mesh(\cU_a)\le\frac{1}{2}L(\wh\cU_{a-1})$,
$\mesh(\wh\cU_a)\le\mesh(\cU_0)$,
$L(\wh\cU_a)\ge\min\{L(\cU_a),\frac{1}{2}L(\wh\cU_{a-1})\}$.
Then using (2), we have
$$\mesh(\cU_{a+1})\le\frac{1}{2}
  \min\{L(\cU_a),\frac{1}{2}L(\wh\cU_{a-1})\}
  \le\frac{1}{2}L(\wh\cU_a).$$
Applying Lemma~\ref{lem:union} to the pair of families
$\wh\cU_a$, $\cU_{a+1}$,
we obtain an open covering
$\wh\cU_{a+1}$
of
$Z_0\cup\dots\cup Z_{a+1}$
with multiplicity
$\le n+1$
and with
$\mesh(\wh\cU_{a+1})\le\max
  \{\mesh(\wh\cU_a),\mesh(\cU_{a+1})\}
  \le\mesh(\cU_0)$
and
$L(\wh\cU_{a+1})\ge\min\{L(\cU_{a+1}),
   \frac{1}{2}L(\wh\cU_a)\}$.

Proceeding by induction and using (3), we obtain an open covering
$\cU=\wh\cU_N$
of
$Z$
of multiplicity
$\le n+1$
with
$\mesh(\cU)\le\mesh(\cU_0)\le\de\tau/2$
and
$L(\cU)\ge\min\{L(\cU_N),\frac{1}{2}L(\cU_{N-1}),\dots,
 \frac{1}{2^N}L(\cU_0)\}\ge(l/\la)\tau$.
Because we can choose
$\tau>0$
arbitrarily small and the constants
$\de$, $\la$, $l$
are independent of
$\tau$,
this shows that
$\cdim Z\le n$.
\qed

\begin{proof}[Proof of Corollary~\ref{cor:selfsim}]
We have
$\dim Z\le\cdim Z$
for every metric space
$Z$.
By Theorem~\ref{thm:main},
$\cdim Z<\infty$
and
$\cdim Z\le\dim Z$,
hence
$\cdim Z=\dim Z$
is finite.
\end{proof}

\subsection{The capacity dimension versus the
topological one}\label{subsect:cdimvdim}

The following proposition allows to construct various
examples of compact metric spaces with the capacity
dimension arbitrarily larger than the topological
dimension.

\begin{pro}\label{pro:locsimcdim} Let
$X$, $Y$
be bounded metric spaces such that for every
$\ep>0$
there is
$A\sub X$
and a homothety
$h_\ep:A\to Y$
is
$\ep$-dense
image,
$\dist(y,h_\ep(A))<\ep$
for every
$y\in Y$.
Then
$\cdim X\ge\dim Y$.
\end{pro}

\begin{proof} We can assume that
$\dim Y>0$,
in particular,
$\diam Y>0$.
Then, we have
$\la(\ep)\ge\la_0>0$
as
$\ep\to 0$
for the coefficient
$\la(\ep)$
of the homothety
$h_\ep$,
because
$X$
is bounded.

Assume that
$n=\cdim X<\dim Y$.
There is
$\de>0$
such that for every sufficiently small
$\tau>0$
there is an open covering
$\cU_\tau$
of
$X$
of multiplicity
$\le n+1$
with
$\mesh(\cU_\tau)\le\tau$
and
$L(\cU_\tau)\ge\de\tau$.

Using the estimate
$\la(\ep)\ge\la_0$,
we can find
$\tau=\tau(\ep)$
such that
$\la(\ep)\tau(\ep)\to 0$
as
$\ep\to 0$
and
$\de\la(\ep)\tau(\ep)\ge 4\ep$.
Then for the covering
$\cV_\ep=h_\ep(\cU_\tau)$
of
$h_\ep(A)$,
we have
$\mesh(\cV_\ep)\le\la(\ep)\tau(\ep)$
and
$L(\cV_\ep)\ge\de\la(\ep)\tau(\ep)$.
Furthermore,
$m(\cV_\ep)\le n+1$.
Therefore, the family
$\cV_{\ep}'=B_{-2\ep}(\cV_\ep)$
still covers
$f(A)$.
Taking the
$\ep$-neighborhood
in
$Y$
of every
$V\in\cV_{\ep}'$,
we obtain an open covering
$\cV$
of
$Y$
with
$\mesh(\cV)\le\mesh(\cV_\ep)\to 0$
as
$\ep\to 0$.

Let us estimate the multiplicity of
$\cV$.
Assume that
$y\in Y$
is a common point of members
$V_j\in\cV$, $j\in J$.
By the definition of
$\cV$,
for every
$j\in J$,
there is
$a_i\in A$
such that
$f(a_j)\in V_j'\in\cV_{\ep}'$
and
$|f(a_j)y|<\ep$.
Then, the mutual distances of the points
$f(a_j)$, $j\in J$,
are
$<2\ep$.
Because
$V_j'=B_{-2\ep}(U_j)$
for
$U_j\in\cV_\ep$,
we see that every point
$f(a_j)$, $j\in J$,
is contained in every
$U_i$, $i\in J$,
and therefore
$|J|\le n+1$
because the multiplicity of
$\cV_\ep$
is at most
$n+1$.
Hence,
$m(\cV)\le n+1$
and
$\dim Y\le n$,
a contradiction.
\end{proof}

As an application, we obtain the following examples.
Let
$Z=\{0\}\cup\set{1/m}{$m\in\N$}$.
Then
$\cdim Z^n=n$
for any
$n\ge 1$,
while
$\dim Z^n=0$.
Indeed, the spaces
$X=Z^n$
and
$Y=[0,1]^n$,
obviously, satisfy the condition of
Proposition~\ref{pro:locsimcdim},
thus
$\cdim Z^n\ge\dim Y=n$
(we have the equality here because
$Z^n\sub Y$).
For
$n=1$,
this example is given in \cite{LS} in context of the
Assouad-Nagata dimension.

Combining with Corollary~\ref{cor:selfsim} and
quasi-symmetry invariance of the capacity dimension,
see \cite{Bu1}, we obtain: The space
$Z^n$
is not quasi-symmetric to any locally self-similar space
for any
$n\ge 1$.

Further examples. Take any monotone sequence of positive
$\ep_k\to 0$, $\ep_1=1/3$,
and repeat the construction of the standard ternary
Cantor set
$K\sub[0,1]$,
only removing at every
$k$-th
step,
$k\ge 1$,
instead of the
$(1/3)^k$-length
segments, the middle segments of length
$s_k=\ep_kl_k$, $l_1=1$,
where the length
$l_{k+1}$
of the segments obtained after processing the
$k$-th
step is defined recurrently by
$2l_{k+1}+s_k=l_k$.
The resulting compact space
$K_a\sub[0,1]$
is homeomorphic to
$K$.
One easily sees that
$\cdim K=0$.
However,
$X=K_a$
and
$Y=[0,1]$
satisfy the condition of Proposition~\ref{pro:locsimcdim},
thus
$\cdim K_a=1$,
while
$\dim K_a=0$.

Similarly, one can construct `exotic' Sierpinski carpets,
Menger curves etc with the capacity dimension strictly
bigger than the topological dimension. Any of those
compact metric spaces is not quasi-symmetric to any
locally self-similar space, in particular, it
is not quasi-symmetric be the boundary (viewed with
a visual metric) of a hyperbolic group, see
Theorem~\ref{thm:cdimhypgroup}. To compare, it is well
known that the boundary at infinity of a typical
hyperbolic group is homeomorphic to the Menger curve.

\section{Proof of Theorem~\ref{thm:main2}}

We actually prove that under the condition of
Theorem~\ref{thm:main2}, the following three
dimensions are finite and coincide
$$\asdim X=\Cdim X=\dim Y.$$

We have
$\asdim X\le\Cdim X$
for every metric space
$X$,
and by Lemmas~\ref{lem:scaledoubling} and
\ref{lem:finitedim},
$\Cdim X$
is finite, because
$X$
is asymptotically similar to the compact space
$Y$.
We first show that
$\dim Y\le\asdim X$.
We let
$\asdim X=N$.
Then for every sufficiently large
$\tau>1$,
there exists an open covering
$\cV$
of
$X$
with multiplicity
$\le N+1$,
$\mesh(\cV)<\infty$
and
$L(\cV)\ge\tau$.

There are constants
$\La_0$,$\la\ge 1$,
such that for every sufficiently large
$R>1$
and every
$V\sub X$
with
$\diam V\le R/\La_0$,
there is a
$\la$-quasi-homothetic
map
$h_V:Y\to X$
with coefficient
$R$,
whose image contains an isometric copy
$V'\sub X$
of
$V$, $V'\sub h_V(Y)$.

Given
$\ep>0$,
we take a sufficiently large
$R>1$
with
$\frac{\la}{R}\mesh(\cV)<\ep$,
satisfying the condition above.
Then, there is a
$\la$-quasi-homothetic
map
$h:Y\to X$
with coefficient
$R$.
The family
$\cU=h^{-1}(\cV)$
is an open covering of
$Y$
with multiplicity
$\le N+1$
and
$\mesh(\cU)\le\frac{\la}{R}\mesh(\cV)<\ep$.
Hence,
$\dim Y\le\asdim X$.

It remains to show that
$\Cdim X\le\dim Y$.
We already know that the topological dimension of
$Y$
and the asymptotic capacity dimension of
$X$
are finite, and we let
$\dim Y=n$, $\Cdim X=N$.
Starting from this point, the proof is completely
parallel to that of Theorem~\ref{thm:main}.

According to the definition of
$\Cdim$,
there is
$\de\in(0,1)$
such that for every sufficiently large
$R$
there exists a
$(N+1)$-colored
open covering
$\cV=\cup_{a\in A}\cV^a$
of
$X$
with
$\mesh(\cV)\le R$
and
$L(\cV)\ge \de R$.
It is convenient to take
$A=\{0,\dots,N\}$
as the color set.

Using that
$Y$
is compact and
$\dim Y=n$,
we find for every
$a\in A$
a finite open covering
$\wt\cU_a$
of
$Y$
with multiplicity
$m(\wt\cU_a)\le n+1$
such the following holds:

\begin{itemize}
\item[(i)] $\mesh(\wt\cU_0)\le\frac{\de}{2\la\La_0}$;
\item[(ii)] $\mesh(\wt\cU_{a+1})
    \le\frac{1}{2\la^2}\min\{L(\wt\cU_a),\mesh(\wt\cU_a)\}$
for every
$a\in A$, $a\le N-1$.
\end{itemize}
Then
$l:=\min\{L(\wt\cU_N),\frac{1}{2}L(\wt\cU_{N-1}),
   \dots,\frac{1}{2^N}L(\wt\cU_0)\}>0$
and
$\mesh(\wt\cU_a)\le\frac{\de}{2\la\La_0}$
for every
$a\in A$.

We put
$t=\de R/2$,
and for every
$V\in\cV$,
consider the smaller subset
$V_t=B_{-t}(V)$.
Then, the sets
$X_a=\cup_{V\in\cV^a}V_t\sub X$, $a\in A$,
cover
$X$, $X=\cup_{a\in A}X_a$,
because
$L(\cV)\ge\de R=2t$.

Given
$V\in\cV$,
there is a
$\la$-quasi-homothetic
map
$h_V:Y\to X$
with coefficient
$\La_0R$
and
$V'\sub h_V(Y)$
for some isometric copy
$V'\sub X$
of
$V$.
Taking the inverse, we obtain a
$\la$-quasi-isometric
map
$f_V:V\to Y$
with coefficient
$(\La_0R)^{-1}$.

Now, for every
$a\in A$, $V\in\cV^a$
consider the family
$\cU_{a,V}=\set{f_V^{-1}(\wt U)}
  {$\wt U\in\wt\cU_a, f_V(V_t)\cap\wt U\neq\es$}$,
which is obviously a covering of
$V_t$
with multiplicity
$\le n+1$.

Note that every
$U\in\cU_{a,V}$, $U=f_V^{-1}(\wt U)$
is contained in
$V$
because
$\dist(v,X\sm V)>t$
for every
$v\in V_t$
and
$\diam U\le\la\La_0R\diam\wt U\le t$.
Thus the family
$\cU_{a,V}$
is contained in
$V$.
Now, the family
$\cU_a=\cup_{V\in\cV^a}\cU_{a,V}$
covers the set
$X_a$
of the color
$a$,
and it has the following properties

\begin{itemize}
\item[(1)] for every
$a\in A$,
the multiplicity of
$\cU_a$
is at most
$n+1$;

\item[(2)] $\mesh(\cU_{a+1})
   \le\frac{1}{2}\min\{L(\cU_a),\mesh(\cU_a)\}$
for every
$a\in A$, $a\le N-1$
($L(\cU_a)$
means the Lebesgue number of
$\cU_a$
as a covering of
$X_a$);

\item[(3)] $\mesh(\cU_a)\le\la\La_0R\mesh(\wt\cU_a)$
and
$L(\cU_a)\ge\La_0R\cdot L(\wt\cU_a)/\la$
for every
$a\in A$.
\end{itemize}

For (1) and (3), the argument is literally the same as
in the proof of Theorem~\ref{thm:main}. Furthermore,
for every
$a\in A$, $a\le N-1$,
and every
$U\in\cU_{a+1}$,
we have
\begin{eqnarray*}
  \diam U&\le&\la\La_0R\mesh(\wt\cU_{a+1})
  \le\frac{\La_0R}{2\la}\min\{L(\wt\cU_a),\mesh(\wt\cU_a)\}\\
  &\le&\frac{1}{2}\min\{L(\cU_a),\mesh(\cU_a)\}
\end{eqnarray*}
by Lemma~\ref{lem:distort}, hence, (2).

Now, we put
$\wh\cU_{-1}=\{X\}$, $\wh\cU_0=\cU_0$
and assume that for some
$a\in A$,
we have already constructed families
$\wh\cU_0,\dots,\wh\cU_a$
so that
$\wh\cU_a$
covers
$X_0\cup\dots\cup X_a$
with multiplicity
$\le n+1$
and
$\mesh(\cU_a)\le\frac{1}{2}L(\wh\cU_{a-1})$,
$\mesh(\wh\cU_a)\le\mesh(\cU_0)$,
$L(\wh\cU_a)\ge\min\{L(\cU_a),\frac{1}{2}L(\wh\cU_{a-1})\}$.
Then using (2), we have
$$\mesh(\cU_{a+1})\le\frac{1}{2}
  \min\{L(\cU_a),\frac{1}{2}L(\wh\cU_{a-1})\}
  \le\frac{1}{2}L(\wh\cU_a).$$
Applying Lemma~\ref{lem:union} to the pair of families
$\wh\cU_a$, $\cU_{a+1}$,
we obtain an open covering
$\wh\cU_{a+1}$
of
$X_0\cup\dots\cup X_{a+1}$
with multiplicity
$\le n+1$
and with
$\mesh(\wh\cU_{a+1})\le\max
  \{\mesh(\wh\cU_a),\mesh(\cU_{a+1})\}
  \le\mesh(\cU_0)$
and
$L(\wh\cU_{a+1})\ge\min\{L(\cU_{a+1}),
   \frac{1}{2}L(\wh\cU_a)\}$.

Proceeding by induction and using (3), we obtain an open covering
$\cU=\wh\cU_N$
of
$X$
of multiplicity
$\le n+1$
with
$\mesh(\cU)\le\mesh(\cU_0)\le\de R/2$
and
$L(\cU)\ge\min\{L(\cU_N),\frac{1}{2}L(\cU_{N-1}),\dots,
 \frac{1}{2^N}L(\cU_0)\}\ge(l\La_0/\la)R$.
Because we can choose
$R$
arbitrarily large and the constants
$\de$, $\la$, $\La_0$, $l$
are independent of
$R$,
this shows that
$\Cdim X\le n$.
\qed

\section{Applications}

\subsection{Capacity dimension of the boundary at infinity
of a hyperbolic group}\label{subsect:cdimbdhypgroup}

Here, we describe a large class of hyperbolic spaces whose
boundary at infinity is locally self-similar and prove
a generalization of Theorem~\ref{thm:cdimhypgroup}.

Recall necessary facts from the hyperbolic spaces theory.
For more details the reader may consult e.g. \cite{BoS}.
We also assume that the reader is familiar with notions
like of a geodesic metric space, a triangle, a geodesic ray etc.

Let
$X$
be a geodesic metric space. We use notation
$xx'$
for a geodesic in
$X$
between
$x$, $x'\in X$,
and
$|xx'|$
for the distance between them. For
$o\in X$
and for
$x$, $x'\in X$,
put
$(x|x')_o=\frac{1}{2}(|xo|+|x'o|-|xx'|)$.
The number
$(x|x')_o$
is nonnegative by the triangle inequality, and it is
called the Gromov product of
$x$, $x'$
w.r.t.
$o$.

\begin{lem}\label{lem:compare} Let
$o$, $g$, $x'$, $x''$
be points of a metric space
$X$
such that
$(x'|g)_o$, $(x''|g)_o\ge |og|-\si$
for some
$\si\ge 0$.
Then
$$(x'|x'')_o\le(x'|x'')_g+|og|\le(x'|x'')_o+2\si.$$
\end{lem}

\begin{proof} The left hand inequality immediately follows
from the triangle inequality: because
$|ox'|\le|og|+|gx'|$
and
$|ox''|\le|og|+|gx''|$,
we have
$(x'|x'')_o\le(x'|x'')_g+|og|$.

Next, we note that
$(x'|o)_g=|og|-(x'|g)_o\le\si$.
This yields
$|x'o|=|og|+|gx'|-2(x'|o)_g\ge|og|+|gx'|-2\si$
and similarly
$|x''o|\ge|og|+|gx''|-2\si$.
Now, the right hand inequality follows.
\end{proof}

A geodesic metric space
$X$
is called
$\de$-{\em hyperbolic},
$\de\ge 0$,
if for any triangle
$xyz\sub X$
the following holds: Let
$y'\in xy$, $z'\in xz$
be points
with
$|xy'|=|xz'|\le(y|z)_x$.
Then
$|y'z'|\le\de$.
In this case, the
$\de$-inequality
$$(x|y)_o\ge\min\{(x|z)_o,(z|y)_o\}-\de$$
holds for every base point
$o\in X$
and all
$x$, $y$, $z\in X$.
A geodesic space is {\em (Gromov) hyperbolic} if it is
$\de$-hyperbolic
for some
$\de\ge 0$.

Let
$X$
be a
$\de$-hyperbolic
space and
$o\in X$
be a base point. A sequence of points
$\{x_i\}\sub X$
{\em converges to infinity,} if
$$\lim_{i,j\to\infty}(x_i|x_j)_o=\infty.$$
Two sequences
$\{x_i\}$, $\{x_i'\}$
that converge to infinity are {\em equivalent} if
$$\lim_{i\to\infty}(x_i|x_i')_o=\infty.$$
Using the
$\de$-inequality,
one easily sees that this defines an equivalence relation
for sequences in
$X$
converging to infinity. The {\em boundary at infinity}
$\di X$
of
$X$
is defined as the set of equivalence classes
of sequences converging to infinity. Every isometry of
$X$
canonically extends to a bijection of
$\di X$
on itself, and we use the same notation for the extension.

Every geodesic ray in
$X$
represents a point at infinity. Conversely, if
a geodesic hyperbolic space
$X$
is proper (i.e. closed balls in
$X$
are compact), then
every point at infinity is represented by a
geodesic ray.

The Gromov product extends to
$X\cup\di X$
as follows. For points
$\xi$, $\xi'\in\di X$,
it is defined by
$$(\xi|\xi')_o=\inf\liminf_{i\to\infty}(x_i|x_i')_o,$$
where the infimum is taken over all sequences
$\{x_i\}\in\xi$, $\{x_i'\}\in\xi'$.
Note that
$(\xi|\xi')_o$
takes values in
$[0,\infty]$,
and that
$(\xi|\xi')_o=\infty$
if and only if
$\xi=\xi'$.

Similarly, the Gromov product
$$(x|\xi)_o=\inf\liminf_{i\to\infty}(x|x_i)_o$$
is defined for any
$x\in X$, $\xi\in\di X$,
where the infimum is taken over all sequences
$\{x_i\}\in\xi$.

Furthermore, for any
$\xi$, $\xi'\in X\cup\di X$
and for arbitrary sequences
$\{x_i\}\in \xi$, $\{x'_i\}\in \xi'$,
we have
$$(\xi|\xi')_o\leq \liminf_{i\to\infty}{(x_i|x'_i)_o}
\leq \limsup_{i\to\infty}{(x_i|x'_i)_o} \leq (\xi|\xi')_o + 2\de.$$
Moreover, the
$\de$-inequality
holds in
$X\cup\di X$,
$$(\xi|\xi'')_o\ge\min\{(\xi|\xi')_o,(\xi'|\xi'')_o\}-\de$$
for every
$\xi$, $\xi'$, $\xi''\in X\cup\di X$.

A metric
$d$
on the boundary at infinity
$\di X$
of
$X$
is said to be {\em visual}, if there are
$o\in X$, $a>1$
and positive constants
$c_1$, $c_2$,
such that
$$c_1a^{-(\xi|\xi')_o}\le d(\xi,\xi')\le c_2a^{-(\xi|\xi')_o}$$
for all
$\xi$, $\xi'\in\di X$.
In this case, we say that
$d$
is the visual metric with respect to the base point
$o$
and the parameter
$a$.
The boundary at infinity is bounded and complete w.r.t.
any visual metric, moreover, if
$X$
is proper then
$\di X$
is compact. If
$a>1$
is sufficiently close to 1, then a visual metric with
respect to
$a$
does exist.

A metric space
$X$
is {\em cobounded} if there is a bounded subset
$A\sub X$
such that the orbit of
$A$
under the isometry group of
$X$
covers
$X$.

\begin{pro}\label{pro:ssboundary} The boundary at infinity
of every cobounded, hyperbolic, proper, geodesic space
$X$
is locally self-similar with respect to any visual metric.
\end{pro}

\begin{proof} We can assume that the geodesic space
$X$
is
$\de$-hyperbolic,
$\de\ge 0$,
and that a visual metric
$d$
on
$\di X$
satisfies
$$c^{-1}a^{-(\xi|\xi')_o}\le d(\xi,\xi')\le ca^{-(\xi|\xi')_o}$$
for some base point
$o\in X$,
some constants
$c\ge 1$, $a>1$
and all
$\xi$, $\xi'\in\di X$.
Note that then
$\diam\di X\le c$.

There is
$\rho>0$
such that the orbit of the ball
$B_\rho(o)$
under the isometry group of
$X$
covers
$X$.
Now, we put
$\la =c^2a^{\rho+4\de}$.
Then
$$\La_0=\min\{1,\diam\di X/\la\}\le 1/c.$$
Fix
$R>1$
and consider
$A\sub\di X$
with
$\diam A\le\La_0/R$.
For each
$\xi$, $\xi'\in A$,
we have
$$(\xi|\xi')_o\ge\log_a\frac{R}{c\La_0}\ge\log_aR.$$

We fix
$\xi\in A$.
Since
$X$
is proper, there is a geodesic ray
$o\xi\sub X$
representing
$\xi$.
We take
$g\in o\xi$
with
$a^{|og|}=R$.
Then using the
$\de$-inequality,
we obtain for every
$\xi'\in A$
$$(\xi'|g)_o\ge\min\{(\xi'|\xi)_o,(\xi|g)_o\}-\de
  =|og|-\de$$
because
$(\xi|g)_o=|og|$.

For arbitrary
$\xi'$, $\xi''\in A$,
consider sequences
$\{x'\}\in\xi'$, $\{x''\}\in\xi''$
such that
$(x'|x'')_g\to(\xi'|\xi'')_g$.
We can assume without loss of generality that
$(x'|g)_o$, $(x''|g)_o\ge|og|-\de$
because possible errors in these estimates disappear
while taking the limit, see below.

Applying Lemma~\ref{lem:compare} to the points
$o$, $g$, $x'$, $x''\in X$
and
$\si=\de$,
we obtain
$$(x'|x'')_o-|og|\le(x'|x'')_g\le(x'|x'')_o-|og|+2\de.$$
Passing to the limit, this yields
$$(\xi'|\xi'')_o-|og|\le(\xi'|\xi'')_g\le(\xi'|\xi'')_o-|og|+4\de.$$
There is an isometry
$f:X\to X$
with
$|of(g)|\le\rho$.
Then
$$(\xi'|\xi'')_g-\rho\le(f(\xi')|f(\xi''))_o
  \le(\xi'|\xi'')_g+\rho$$
because the Gromov products with respect to different
points differ one from each other at most by the
distance between the points. The last two double inequalities
give
$$(\xi'|\xi'')_o-|og|-\rho\le(f(\xi')|f(\xi''))_o
  \le(\xi'|\xi'')_o-|og|+\rho+4\de,$$
and therefore,
$$c^{-2}a^{-(\rho+4\de)}Rd(\xi',\xi'')\le
  d(f(\xi'),f(\xi''))\le c^2a^{\rho}Rd(\xi',\xi'').$$
This shows that
$f:A\to\di X$
is
$\la$-quasi-homothetic
with coefficient
$R$
and hence
$\di X$
is locally self-similar.
\end{proof}

Now, Corollary~\ref{cor:selfsim} and
Proposition~\ref{pro:ssboundary} give the following

\begin{thm}\label{thm:cdimhypspace} The capacity
dimension of the boundary at infinity of every cobounded,
hyperbolic, proper, geodesic space
$X$
coincides with the topological dimension,
$\cdim\di X=\dim\di X$.
\qed
\end{thm}

The class of spaces satisfying the condition of
Theorem~\ref{thm:cdimhypspace} is very large.
It includes in particular all symmetric rank one spaces of
noncompact type (i.e. the real, complex, quaternionic hyperbolic
spaces and the Cayley hyperbolic plane), all cocompact Hadamard
manifolds of negative sectional curvature, various hyperbolic
buildings etc. The most important among them is the class
of (Gromov) hyperbolic groups, and Theorem~\ref{thm:cdimhypgroup}
is a particular case of Theorem~\ref{thm:cdimhypspace}.

Note that for many such spaces, the boundary at infinity
is fractal in the sense that its Hausdorff dimension with
respect of a natural visual metric is larger than the
topological dimension. This is true e.g. for the complex,
quaternionic hyperbolic spaces, the Cayley hyperbolic plane,
for the Fuchsian hyperbolic buildings, see \cite{Bou},
and for the hyperbolic graph surfaces, see \cite{Bu3}.

\subsection{The asymptotic dimension of a hyperbolic group}
\label{subsect:asdimhypgroup}

Theorem~\ref{thm:cdimhypspace} is the decisive step in
the proof of the following

\begin{thm}\label{thm:asdim} The asymptotic dimension of
every cobounded, hyperbolic, proper, geodesic space
$X$
equals topological dimension of its boundary at infinity plus 1,
$$\asdim X=\dim\di X+1.$$
\end{thm}

Theorem~\ref{thm:gromovconjecture} is a particular case
of Theorem~\ref{thm:asdim}. The fact that
$\asdim X$
as well as
$\dim\di X$
are finite,
$\asdim X$, $\dim\di X<\infty$,
is well known, it follows e.g. from \cite{BoS} (for hyperbolic
groups, there is an alternative proof \cite[Theorem~9.25]{Ro}).
Our contribution is that we prove the optimal estimate,
$\asdim X\le\dim\di X+1$.

The estimate from below,
$$\asdim X\ge\dim\di X+1,$$
or at least the idea of its proof is also well known, see
\cite[1.$\text{E}_1'$]{Gr}. More stronger estimates
of different types are obtained in \cite{Sw} and \cite{BS2}
respectively. For convenience of the reader, we give
a simplified version of arguments from \cite{BS2}
adapted to the asymptotic dimension.

Let
$Z$
be a bounded metric space. Assuming that
$\diam Z>0$,
we put
$\mu=\pi/\diam Z$
and note that
$\mu|zz'|\in[0,\pi]$
for every
$z$, $z'\in Z$.
Recall that the hyperbolic cone
$\cone(Z)$
over
$Z$
is the space
$Z\times[0,\infty)/Z\times\{0\}$
with metric defined as follows. Given
$x=(z,t)$, $x'=(z',t')\in\cone(Z)$
we consider a triangle
$\ov o\,\ov x\,\ov x'\sub\hyp^2$
with
$|\ov o\,\ov x|=t$, $|\ov o\,\ov x'|=t'$
and the angle
$\angle_{\ov o}(\ov x,\ov x')=\mu|zz'|$.
Now, we put
$|xx'|:=|\ov x\,\ov x'|$.
In the degenerate case
$Z=\{\pt\}$,
we define
$\cone(Z)=\{\pt\}\times[0,\infty)$
as the metric product. The point
$o=Z\times\{0\}\in\cone(Z)$
is called the {\em vertex} of
$\cone(Z)$.

\begin{pro}\label{pro:asdimbelow} For every geodesic hyperbolic
space
$X$
we have
$$\asdim X\ge\dim\di X+1.$$
\end{pro}

\begin{proof} The same argument as in \cite[Proposition~6.2]{Bu1}
shows that
the hyperbolic cone
$\cone(Z)$
over
$Z=\di X$,
taken with some visual metric, can be quasi-isometrically
(actually, roughly similarly) embedded in
$X$
because
$X$
is geodesic. Thus
$\asdim X\ge\asdim\cone(Z)$,
and we show that
$\asdim\cone(Z)\ge\dim Z+1$.

The {\em annulus}
$\an(Z)\sub\cone(Z)$
consists of all
$x\in\cone(Z)$
with
$1\le|xo|\le 2$.
Clearly,
$\an(Z)$
is homeomorphic to
$Z\times I$, $I=[0,1]$.
According to a well known result from the dimension
theory (see \cite{Al}), the topological dimension
\[\dim\an(Z)=\dim Z+1.\]

Consider the sequence of contracting homeomorphisms
$F_k:\cone(Z)\to\cone(Z)$
given by
$F_k(z,t)=(z,\frac{1}{k}t)$, $(z,t)\in\cone(Z)$, $k\in\N$.
Given a uniformly bounded covering
$\cU$
of
$\cone(Z)$,
the coverings
$\cU_k=F_k(\cU)\cap\an(Z)$
of the annulus
$\an(Z)$
have arbitrarily small
$\mesh$
as
$k\to\infty$.
Therefore,
$\asdim\cone(Z)\ge\dim\an(Z)$,
and the estimate follows.
\end{proof}

\begin{proof}[Proof of Theorem~\ref{thm:asdim}]
The estimate from below
$$\asdim X\ge\dim\di X+1$$
follows from Proposition~\ref{pro:asdimbelow}. By the main
result of \cite{Bu1}, see also \cite{Bu2} for another proof,
we have
$\asdim X\le\cdim\di X+1$
(the space
$X$
is certainly visual, i.e.
$X$
satisfies the condition of the cited theorems). Now,
the estimate from above,
$$\asdim X\le\dim\di X+1,$$
follows from Theorem~\ref{thm:cdimhypspace}.
\end{proof}

\subsection{Embedding and nonembedding results}
\label{subsect:embnonemb}

Combining Theorem~\ref{thm:cdimhypspace} with the main
result of \cite{Bu2}, we obtain

\begin{thm}\label{thm:embdspace} Every cobounded,
hyperbolic, proper, geodesic space
$X$
admits a quasi-isometric embedding
$X\to T_1\times\dots\times T_n$
into the
$n$-fold
product of simplicial metric trees
$T_1,\dots,T_n$
with
$n=\dim\di X+1$.
\qed
\end{thm}

\noindent
(Every space
$X$
satisfying the condition of Theorem~\ref{thm:embdspace}
is certainly visual, i.e.
$X$
also satisfies the condition of \cite[Theorem~1.1]{Bu2}).
For example, the complex hyperbolic plane
$\C\hyp^2$
admits a quasi-isometric embedding into the
$4$-fold
product of simplicial metric trees etc.

Theorem~\ref{thm:embd} is a particular case of
Theorem~\ref{thm:embdspace}.

Theorem~\ref{thm:cdimhypspace} has applications also to
nonembedding results. For example, let
$X^n$
be a universal covering of a compact Riemannian
$n$-dimensional,
$n\ge 2$
manifold with nonempty geodesic boundary and constant
sectional curvature
$-1$.
Then
$X^n$
satisfies the condition of Theorem~\ref{thm:cdimhypspace},
and hence
$\dim\di X^n=\cdim\di X^n$.
Note that
$X^n$
can be obtained from the real hyperbolic space
$\hyp^n$
by removing a countable collection of disjoint open
half-spaces, and
$\di X^n\sub S^{n-1}$
is a compact nowhere dense subset obtained by removing
a countable collection of disjoint open balls.
In particular, for
$n=2$,
$\di X^n\sub S^1$
is a Cantor set, for
$n=3$,
$\di X^n\sub S^2$
is a Sierpinski carpet, and for
$n\ge 4$,
$\di X^n\sub S^{n-1}$
is a higher dimensional version of a Sierpinski carpet.
Thus
$\dim\di X^n=n-2$.

The space
$X^n$
contains isometrically embedded copies of
$\hyp^{n-1}$
as the boundary components, thus by \cite{BF}, the
$k$-fold
product
$X^n\times\dots\times X^n$,
$k\ge 1$
factors, contains quasi-isometrically embedded
$\hyp^p$
for
$p=k(n-2)+1$.

\begin{thm}\label{thm:nonembd} Let
$X^n$
be the space as above,
$Y_k^n=X^n\times\dots\times X^n$
be the
$k$-fold
product,
$k\ge 1$.
Then there is no quasi-isometric embedding
$$\hyp^p\to Y_k^n\times\R^m$$
for
$p>k(n-1)$
and any
$m\ge 0$.
\end{thm}

For example, the Theorem says that there is no way to
embed quasi-isometrically
$\hyp^5$
into
$X^3\times X^3\times\R^m$
for any
$m\ge 0$
though there is a quasi-isometric
$\hyp^3\to X^3\times X^3$.
In general, for an arbitrary
$n\ge 2$,
there is a quasi-isometric
$\hyp^p\to X^n\times X^n$
for
$p=2n-3$
and there is no quasi-isometric
$\hyp^p\to X^n\times X^n\times\R^m$
for
$p=2n-1$
and an arbitrary
$m\ge 0$.

In the case
$n=2$,
the space
$X^2$
is quasi-isometric to the binary tree
$T$
whose edges all have length 1 because
$X^2$
covers a compact hyperbolic surface with nonempty
geodesic boundary. By \cite{DS}, there is a quasi-isometric
embedding
$\hyp^2\to T\times T$,
and hence there is a quasi-isometric
$\hyp^p\to X^n\times X^n$
in the remaining case
$p=2n-2$
if
$n=2$.
For
$n\ge 3$,
the question whether there is a quasi-isometric
$\hyp^{2n-2}\to X^n\times X^n$
remains open. Moreover, the same question is open
for quasi-isometric
$$\hyp^{k(n-1)}\to Y_k^n,\ n,k\ge 3.$$

Certainly, there is a huge range of possibilities in
variation of this theme, e.g. taking as the target
space the product with different dimensions of factors,
considering other than
$X$
spaces, replacing
$\hyp^p$
as source space etc. However, everything what is known here
on nonembedding side is covered by ideas of the proof of
Theorem~\ref{thm:nonembd}.

\begin{proof}[Proof of Theorem~\ref{thm:nonembd}] The key ingredient
of the proof is the notion of the hyperbolic dimension of a metric
space
$X$,
$\hypdim X$,
which is introduced in \cite{BS2} and has the following properties
\begin{itemize}
\item[(1)] monotonicity: let
$f:X\to X'$
be a quasi-isometric map between metric spaces
$X$, $X'$.
Then
$\hypdim X\le\hypdim X'$;
\item[(2)] the product theorem: for any metric spaces
$X_1$, $X_2$,
we have
$$\hypdim(X_1\times X_2)\le\hypdim X_1+\hypdim X_2;$$
\item[(3)] $\hypdim X\le\asdim X$
for every metric space
$X$;
\item[(4)] $\hypdim\R^m=0$
for every
$m\ge 0$.
\end{itemize}

Recall that a metric space
$X$
has {\em bounded growth at some scale}, if for some
constants
$r$, $R$ with
$R>r>0$,
and
$N\in\N$
every ball of radius
$R$
in
$X$
can be covered by
$N$
balls of radius
$r$,
see \cite{BoS}.

According to the main result of \cite{BS2}, for
every geodesic, Gromov hyperbolic space, which has bounded
growth at some scale and whose boundary at infinity
$\di X$
is infinite, one holds
$$\hypdim X\ge\dim\di X+1.$$
In particular,
$\hypdim\hyp^p\ge p$
(actually, the equality here is true). Now, assume that
there is a quasi-isometric embedding
$\hyp^p\to Y_k^n\times\R^m$.
Then by properties (1) -- (4) above, we have
$$p\le\hypdim\hyp^p\le k\cdot\hypdim X^n
  \le k\cdot\asdim X^n.$$
Using Theorem~\ref{thm:embdspace} and properties of
the asymptotic dimension, or using that
$\asdim X^n\le\cdim\di X^n+1$
(\cite{Bu1}) and
Theorem~\ref{thm:cdimhypspace}, we obtain
$\asdim X^n\le\dim\di X^n+1=n-1$.
Hence, the claim.
\end{proof}

\subsection{Examples of strict inequality in the product
theorem for the capacity dimension}

An application of Corollary~\ref{cor:selfsim} is that the strict
inequality in the product theorem for the capacity dimension
holds for some compact metric spaces. With each
$n\in\N$,
one associates a Pontryagin surface
$\Pi_n$
which is a 2-dimensional compactum,
$\dim\Pi_n=2$,
and for coprime
$m$, $n\in\N$
one holds
$\dim(\Pi_m\times\Pi_n)=3$,
that is
$$\dim(\Pi_m\times\Pi_n)<\dim\Pi_m+\dim\Pi_n.$$

According \cite{Dr1}, \cite[Corollary~2.3]{Dr2}, for every prime
$p$,
there is a hyperbolic Coxeter group with a Pontryagin surface
$\Pi_p$
as the boundary at infinity. Taken with any visual metric,
$\Pi_p$
is a locally self-similar space and hence
$\cdim\Pi_p=\dim\Pi_p=2$
by Theorem~\ref{thm:cdimhypgroup}. Obviously, the product
of locally self-similar spaces is locally self-similar, and we
obtain
$$\cdim(\Pi_p\times\Pi_q)<\cdim\Pi_p+\cdim\Pi_q$$
for prime
$p\neq q$.

Unfortunately, the construction in \cite{Dr2} of an appropriate
hyperbolic Coxeter group is implicit.
Here, we give explicit constructions self-similar Pontryagin
surfaces.

\section{Self-similar Pontryagin surfaces}
\label{sect:sspont}

For the notion of the cohomological dimension of a
topological space
$X$
with respect to an abelian group
$G$, $\dim_GX$,
we refer e.g. to the survey \cite{Dr3}. Let
$p$
be a prime number.
A {\em Pontryagin surface
$\Pi_p$}
is a 2-dimensional compact space with
$\dim_{\Q}\Pi_p=\dim_{\Z_q}\Pi_p=1$
for every prime
$q\neq p$
and
$\dim_{\Z_p}\Pi_p=2$.

\begin{thm}\label{thm:sspont} For every prime
$p$,
there exists a Pontryagin surface
$\Pi_p$
with locally self-similar metric.
\end{thm}

Our objective is the existence of a (locally) self-similar
metric space
$\Pi_p$.
For a simple argument, which proves the required cohomological
properties of the compactum
$\Pi_p$
we construct, we refer to \cite[Example~1.9]{Dr3}.

\subsection{Construction}\label{subsect:construct}

By a square we mean a topological space homeomorphic to
$[0,1]\times[0,1]$.
Given a natural
$m\ge 2$,
we consider the
$m$-{\em band}
$B_m$,
which is a 2-dimensional square complex, constructed
as follows. The union
$\wt B_m$
of
$m$
squares with a common side can also be described as
$T_m\times[0,1]$,
where
$T_m$
is the union of
$m$
copies of the segment
$[0,1]$
attached to each other along the common vertex
$0$.
We fix a cyclic permutation
$\tau$
of the segments
$\si\sub T_m$
and define
$B_m=\wt B_m/\set{\si\times 0\equiv\tau(\si)\times 1}{$\si\sub T_m$}$.
In the case
$m=2$,
this gives the usual M\"obius band.

The square complex
$B_m$
consists of
$m$
squares and its boundary
$\d B_m$
as well as its singular locus
$sB_m$
corresponding to the common side of the squares
both are homeomorphic to
$S^1$.

For every
$a$, $b>0$,
there is a well defined intrinsic metric on
$B_m$
with respect to which every square of
$B_m$
is isometric to the Euclidean rectangle whose sides
have length
$a$
and
$b$.
We assume that
$a$
is the length of the common side of the squares.
Then
$\d B_m\sub B_m$
is a geodesic of length
$ma$
and
$sB_m\sub B_m$
is a geodesic of length
$a$.
We use notation
$B_m(a,b)$
for
$B_m$
endowed with this metric. Note that
$B_m(a,b)$
has nonpositive curvature (in Alexandrov sense)
for every
$a$, $b>0$.

\begin{lem}\label{lem:bandpro} Given
$a>0$, $b\ge a\sqrt2/2$,
let
$A=B_m\left(\frac{4a}{m},b\right)$,
$B=[0,a]\times[0,a]$.
Then, there is a 1-Lipschitz map
$$q:(A,\d A)\to(B,\d B)$$
whose restriction to the boundary,
$q|\d A$,
preserves the length of every arc.
\end{lem}

\begin{proof} The boundary
$\d A$
is a closed geodesic in
$A$
of length
$4a$.
We subdivide
$\d A$
into four segments of length
$a$,
called the {\em sides}, and map
$\d A$
onto
$\d B$
the sides and arc length preserving. This defines
$q:\d A\to\d B$,
which is obviously 1-Lipschitz.

Next, we extend
$q$
to the singular locus
$sB_m$
by collapsing it to the middle of the square
$B$.
This extension is still 1-Lipschitz because the distance
between any points
$x\in\d A$
and
$x'\in sB_m$
is at least
$b$
and
$b\ge a\sqrt 2/2$,
which is the maximal distance between the middle of
$B$
and points of
$\d B$.

Finally, we extend already defined
$q$
to
$A$
as the affine map on every radial segment
$xx'\sub A$
with
$x\in\d A$, $x'\in sB_m$
be the (unique) closest to
$x$
point. When
$x$
runs over
$\d A$,
the segments
$xx'$
cover
$A$
so that for different
$x$, $y\in\d A$
the segments
$xx'$
and
$yy'$
have no common interior point. Because
$|xx'|=b$,
the restriction of
$q$
to
$xx'$
is 1-Lipschitz.

So defined
$q$
is smooth outside of the union of
$\d A$, $sB_m$
and the four radial segments corresponding to
the end points of the sides of
$\d A$.
One easily sees that
$|dq|\le 1$
there. Since
$A$
is geodesic, it follows that
$q$
is 1-Lipschitz.
\end{proof}

\subsubsection{First template}
Fix a natural
$m\ge 2$
and an odd
$k$, $k=2l+1$
with
$l\ge 1$.
We define a 2-dimensional square complex
$P=P_{m,k}$
obtained from the square
$Q_k=[0,k]\times[0,k]$
by removing the open middle square
$$q_k=(l,l+1)\times(l,l+1)\sub Q_k$$
and attaching instead the
$m$-band
$B_m$
along a homeomorphism
$\d B_m\to\d q_k$.
We consider the following square structure on
$P$.
The remainder
$Q_k\sm q_k$
consists of
$k^2-1$
squares each of which we subdivide into
$m^2$
subsquares. Furthermore, the
$m$-band
consists of
$m$
squares each of which we represent as the rectangle
$[0,4]\times[0,m]$
with the natural square structure consisting of
$4m$
squares. We assume that the side
$[0,4]\times\{0\}$
corresponds to the singular locus
$sB_m$
of
$B_m$.
Assuming that the gluing homeomorphism
$\d B_m\to\d q_k$
preserves the induced subdivisions of each circle into
$4m$
segments, we obtain the desired square complex structure on
$P$
consisting of
$$s_{m,k}=(k^2-1)m^2+4m^2=(km)^2+3m^2$$
squares. Speaking about squares of
$P$
we mean squares of this square complex structure.

We consider the canonical intrinsic metric on
$P$
with respect to which every square of
$P$
is isometric to the Euclidean square with the
side length
$1/(km)$.
The space
$P$
is nonpositively curved, and the subcomplex
$B_m\sub P$
is convex, isometric to
$B_m(\frac{4}{km},\frac{1}{k})$,
and its boundary is geodesic in
$P$.

The boundary
$\d P$
consists of four sides of length 1 and one can consider
$P$
as the unit square
$[0,1]\times[0,1]$
with the appropriate middle subsquare replaced by
$B_m\left(\frac{4}{km},\frac{1}{k}\right)$.
We have
$b>a\sqrt 2/2$
for
$a=b=1/k$.
Thus applying Lemma~\ref{lem:bandpro}, we obtain
a 1-Lipschitz map
$$q_0^1:P\to [0,1]\times[0,1]$$
which is identical outside of the interior of
$B_m\sub P$.

\subsubsection{Constructing a sequence of polyhedra
$\{P^i\}$}
We construct a sequence of polyhedra
$\{P^i=P_{m,k}^i\}$, $i\ge 1$,
in a way that every
$P^i$
serves as a building block for
$P^{i+1}$
and as such it is called the
$i$-th
{\em template.} Furthermore, every
$P^{i+1}$
consists of one and the same number of building blocks
$P^i$
independent of
$i$.
Every polyhedron
$P^i$
possesses a canonical square complex structure and being
endowed with an intrinsic metric, for which every square
is isometric to a fixed Euclidean square, it is nonpositively
curved in the Alexandrov sense. From ideological side, the
construction is similar to well known constructions of
self-similar compact metric spaces via a family of homotheties.

As the 0-th template we take the unit square,
$P^0=[0,1]\times[0,1]$.
The first template
$P^1=P$
is already described above. It can also be described
as follows. The polyhedron
$P^1$
consists of
$s_{m,k}$
blocks each of which is a
$1/mk$-homothetic
copy of
$P^0$
attached along the boundary to the 1-skeleton
$S$
of
$P^1$.

The polyhedron
$P^2$
is obtained out of
$P^1$
as follows. We remove every of
$s_{m,k}$
open square of
$P^1$,
obtaining again the 1-dimensional complex
$S$,
and replace every removed open square by a
$1/mk$-homothetic
copy of
$P^1$
attaching it along the boundary to
$S$.

Assume that a square polyhedron
$P^i$
is already constructed for
$i\ge 1$.
By assumption, it is considered with the canonical
intrinsic nonpositively curved metric in which every square
is isometric to the Euclidean square of side length
$1/(mk)^i$.
The polyhedron
$P^i$
consists of
$s_{m,k}$
pairwise isometric blocks each of which is
$mk$-homothetic
to the template
$P^{i-1}$.
The boundary
$\d P^i$
is subdivided into four sides, consisting each of
$(mk)^i$
segments of the square structure, and has length 4.

We construct the square polyhedron
$P^{i+1}$,
replacing every of
$s_{m,k}$
open square of
$P^1$
by a
$1/mk$-homothetic
copy of the template
$P^i$,
attaching it to
$S$,
so that they all together form a square complex
structure of
$P^{i+1}$
and define the canonical intrinsic nonpositively curved metric
in which every square is isometric to the Euclidean square
of side length
$1/(mk)^{i+1}$.

Our construction has the following property: for every integer
$j$, $1\le j\le i$,
the polyhedron
$P^{i+1}$
consists of
$s_{m,k}^{i+1-j}$
subblocks each of which is homothetic to the
$j$-th
template
$P^j$
with coefficient
$(mk)^{i+1-j}$.

\begin{lem}\label{lem:diamtemplate} The diameter of
$P^i$, $i\ge 0$,
is bounded above by a constant independent of
$i$,
$$\diam P^i\le d$$
where one can take
$d=d(m,k)=\frac{2m(l+1)}{mk-1}+2$,
(recall
$k=2l+1$).
\end{lem}

\begin{proof} Because the length of the boundary
$\d P^i$
equals 4 for every
$i\ge 1$,
it suffices to estimate
$\de_i=\max\set{\dist(x,\d P^i)}{$x\in P^i$}$
from above independent of
$i$.

For
$i=1$,
the most remote points from the boundary are sitting in
the singular locus of the subcomplex
$B_m\sub P^1$.
Moving along the 1-skeleton
$S$
of
$P^1$,
we find that
$\de_1=m/mk+ml/mk=(l+1)/k$.

The grid
$S$
serves as a skeleton for attaching the blocks while
constructing every
$P^i$
and thus it is isometrically (in the sense of the induced
intrinsic metric) embedded in
$P^i$
for every
$i\ge 1$.
So, to estimate
$\de_i$
we can use paths in
$S$,
namely,
$\dist_S(x,S_0)\le\de_1$
for every
$x\in S$,
where
$S_0\sub S$
is identified with the boundary
$\d P^i$
for every
$i\ge 1$,
and the distance is taken with respect to the intrinsic metric of
$S$.

For
$i=2$,
clearly
$$\de_2\le\de_1/mk+\max_{x\in S}\dist_S(x,S_0)
  \le\de_1/mk+\de_1$$
because
$P^2$
consists of blocks
$1/mk$-homothetic to
$P^1$,
whose boundaries are subsets of
$S$.
Similarly, we recurrently obtain the estimate
$$\de_{i+1}\le\de_i/mk+\de_1
  \le\de_1\sum_{j=0}^i1/(mk)^j,$$
hence, the claim.
\end{proof}

\subsection{The inverse sequence
$\{P^i;q_i^{i+1}\}$}

The bonding map
$q_0^1:P^1\to P^0$
is already described above. By induction,
we obtain the bonding map
$q_i^{i+1}:P^{i+1}\to P^i$
for every
$i\ge 1$
by putting together the maps
$q_{i-1}^i$
defined on the blocks of
$P^{i+1}$.

This map is 1-Lipschitz and it is compatible
with self-similar structure of complexes,
i.e. its restriction to every subblock
$P^j\sub P^{i+1}$, $1\le j\le i$,
coincides with
$q_{j-1}^j$.

The product
$\cP=\prod_{i\ge 0}P^i$
is the set of all sequences
$\set{x_i\in P^i}{$i\ge 0$}$.
The limit of the inverse sequence
$\{P^i;q_i^{i+1}\}$,
$$\Pi=\Pi_{m,k}=\lim_{\longleftarrow}(P^i;q_i^{i+1}),$$
is the subset of
$\cP$
consisting of all sequences
$\{x_i\}$
with
$x_i=q_i^{i+1}(x_{i+1})$
for every
$i\ge 0$.
For every
$j\ge 0$,
we have the projection
$q_j^{\infty}:\Pi\to P^j$
defined by
$q_j^{\infty}(\{x_i\})=x_j$
for every
$\{x_i\}\in\Pi$.
Clearly,
$q_j^{j+1}\circ q_{j+1}^\infty=q_j^\infty$
for every
$j\ge 0$.

The space
$\cP$
is compact in the product topology as the product
of compact spaces, and
$\Pi$
is closed in
$\cP$
because all bonding maps are continuous. Thus
$\Pi$
is compact in the induced topology, which we call
the {\em product topology} of
$\Pi$.
By the definition of the product topology, the map
$q_j^\infty$
is continuous for every
$j\ge 0$
and the product topology is the roughest one among all topologies
with this property.

For each
$\xi=\{x_i\}$, $\xi'=\{x_i'\}\in\Pi$,
the sequence of distances
$|x_ix_i'|$
is bounded by Lemma~\ref{lem:diamtemplate} and
nondecreasing because every bonding map is 1-Lipschitz.
Now, we define a metric on
$\Pi$
by
$$|\xi\xi'|=\lim_i|x_ix_i'|$$
(this metric is of course by no means nonpositively curved
despite the fact that all
$P^i$
are nonpositively curved). The corresponding metric topology
on
$\Pi$
we call the {\em metric topology} of
$\Pi$.

\begin{lem}\label{lem:topology} The metric topology of
$\Pi$
coincides with the product topology.
\end{lem}

\begin{proof} If
$|\xi\xi'|<r$
for some
$\xi=\{x_i\}$, $\xi'=\{x_i'\}\in\Pi$,
then
$|x_ix_i'|\le|\xi\xi'|<r$
for every
$i\ge 0$.
It follows that the projection
$q_i^\infty$
is continuous in the metric topology for every
$i\ge 0$
and thus every open set in the product topology is
open in the metric topology.

Fix
$\xi\in\Pi$, $r>0$
and consider the (open) ball
$B_r(\xi)\sub\Pi$.
We show that there is an open in the product topology
subset which is contained in
$B_r(\xi)$
and contains
$\xi$.
There is
$j\in\N$
such that
$2d/(mk)^j<r$,
where
$d=d(m,k)$
is the upper bound for the diameter of
$P^i$,
see Lemma~\ref{lem:diamtemplate}, and hence for
the diameter of
$\Pi$.

We consider
$x_j=q_j^\infty(\xi)\in P^j$
and take the union
$\ov A$
of all squares of
$P^j$
containing
$x_j$.
Recall that the side length of every such square is
$1/(mk)^j$.
Then, the point
$x_j$
is contained in the interior
$A$
of
$\ov A$,
which open in
$P^j$.
Thus
$B=(q_j^\infty)^{-1}(A)$
is open in the product topology, and
$\xi\in B$.
On the other hand,
$\diam B\le 2\diam C$,
where
$C$
is preimage under
$q_j^\infty$
of a square of
$P^j$.
We have
$\diam C=\diam\Pi/(mk)^j\le d/(mk)^j$
and hence
$\diam B\le 2d/(mk)^j<r$.
It follows that
$B\sub B_r(\xi)$.
Therefore, every open in the metric topology
subset in
$\Pi$
is also open in the product topology.
\end{proof}

\subsection{Self-similarity of the space
$\Pi$}

Recall some notions from the theory of self-similar
metric spaces, see e.g. \cite{Fa}, \cite{Hu}.
A compact metric space
$K$
is said to be {\em self-similar} if there is a finite
collection of homotheties
$f_a:K\to K$, $a\in A$,
with coefficients
$h_a\in(0,1)$
such that
$K\sub\cup_{a\in A}f_a(K)$.

In this case, there is a unique number
$\mu\ge 0$
with
$$\sum_{a\in A}h_a^\mu=1.$$
The number
$\mu$
is called the similarity dimension of
$K$, $\mu=\dim_sK$
(more precisely,
$\mu$
is the similarity dimension of the family
$\{f_a\}_{a\in A}$).
For example, if
$h_a=h$
for all
$a\in A$,
then we have
$$\dim_sK=\frac{\log|A|}{\log(1/h)}.$$
One always has
$\dim_HK\le\dim_sK$
for the Hausdorff dimension
$\dim_HK$
of
$K$.
The collection of the homotheties
$\set{f_a}{$a\in A$}$
satisfies the {\em OSC} (Open Set Condition), if
there is an open set
$U\sub K$
such that
$f_a(U)\sub U$
for all
$a\in A$
and
$f_a(U)\cap f_{a'}(U)=\es$
for all distinct
$a$, $a'\in A$.
In this case, the Hausdorff dimension of
$K$
coincides with the similarity dimension,
$\dim_HK=\dim_sK$
(see e.g. \cite{EG}).

\begin{pro}\label{pro:sspont} The metric space
$\Pi=\Pi_{m,k}$
is compact and self-similar for every integer
$m\ge 2$
and odd
$k=2l+1$.
Furthermore, the corresponding collection
$\set{f_a}{$a\in A$}$
consists of
$|A|=s_{m,k}$
homotheties with coefficients
$h_a=1/(mk)$,
and satisfies the OSC. In particular, the Hausdorff dimension
$$\dim_H\Pi=2+\frac{\log(1+3/k^2)}{\log(mk)}.$$
\end{pro}

\begin{proof} It follows from Lemma~\ref{lem:topology}
that the metric space
$\Pi$
is compact. Recall that the square polyhedron
$P^{i+1}$
consists of
$s_{m,k}=(mk)^2+3m^2$
blocks with disjoint interiors, each of which is
$mk$-homothetic to
$P^i$
for every
$i\ge 0$.

We label the blocks of
$P^i$
by a finite set
$A$, $|A|=s_{m,k}$,
in a way independent of
$i$,
that is compatible with the bonding maps, and fix for each
$i\ge 0$, $a\in A$
a homothety
$f_a^i:P^i\to P^{i+1}$
with coefficient
$h_a=1/(mk)$
whose image is the corresponding block of
$P^{i+1}$.
We can also assume that
\[f_a^i\circ q_i^{i+1}=q_{i+1}^{i+2}\circ f_a^{i+1}\tag{$\ast$}\]
for each
$i\ge 0$, $a\in A$.
Then
$P^{i+1}=\cup_{a\in A}f_a^i(P^i)$
and for the interior
$U^i$
of
$P^i$,
we have
$f_a^i(U^i)\sub U^{i+1}$
while the open sets
$f_a^i(U^i)$, $f_{a'}^i(U^i)$
are disjoint for different
$a$, $a'\in A$.

The equality
$(\ast)$
allows to pass to the limit as
$i\to\infty$,
which yields the collection of homotheties
$f_a:\Pi\to\Pi$, $a\in A$,
with coefficients
$h_a=1/(mk)$
for every
$a\in A$
with the required properties.
\end{proof}

\begin{proof}[Proof of Theorem~\ref{thm:sspont}]
We fix an odd
$k=2l+1\ge 3$
and consider the metric space
$\Pi_p=\Pi_{p,k}$,
which is self-similar by Proposition~\ref{pro:sspont}.
It is an easy exercise to check that every compact self-similar
space is locally self-similar. We have
$2<\dim_H\Pi_p<3$
for the Hausdorff dimension by Proposition~\ref{pro:sspont},
thus
$\dim\Pi_p\le 2$.
Applying the argument from \cite[Example~1.9]{Dr3}, we obtain
$\dim_{\Q}\Pi_p=\dim_{\Z_q}\Pi_p=1$
for every prime
$q\neq p$
and
$\dim_{\Z_p}\Pi_p=2$.
The last equality together with the estimate
$\dim\Pi_p\le 2$
implies that
$\dim\Pi_p=2$.
\end{proof}

\subsection{Asymptotically self-similar Pontryagin surfaces}
\label{subsect:asympont}

We fix a compact self-similar Pontryagin surface
$\Pi=\Pi_{m,k}$
for some integer
$m\ge 2$,
odd
$k=2l+1\ge 3$,
and define a metric space
$\wh\Pi$
as follows. Recall
$\Pi$
consists of
$s_{m,k}=(mk)^2+3m^2$
blocks, each of which is
$(1/mk)$-homothetic
to
$\Pi$
and attached along the 1-skeleton
$S$
of the polyhedron
$P^1$,
and
$S$
is isometrically (in the sense of the induced
intrinsic metric) embedded in
$\Pi$.
For every of
$4m^2$
middle blocks, which are projected into the
$m$-band $B_m\sub P^1$,
their distance to the boundary
$\d\Pi$
is at least
$ml$,
because the projection
$\Pi\to P^1$
is 1-Lipschitz, and the distance from
$B_m$
to the boundary
$\d P^1$
equals
$ml$.

We put
$\la=mk$, $X_0=\Pi$,
and consider the
$\la$-homothetic
copy of
$X_0$,
$X_1=\la X_0$.
The space
$X_1$
consists of
$s_{m,k}$
blocks isometric to
$X_0$,
and we can consider
$X_0$
as a subspace of
$X_1$, $X_0\sub X_1$,
identifying it with some block. Moreover, we take
this block in the middle of
$X_1$
so that
$\dist(X_0,\d X_1)\ge\la ml$.
Taking
$X_i=\la^iX_0$,
we obtain an increasing sequence
$$X_0\sub X_1\sub\dots\sub X_i$$
of metric spaces via appropriate identifications
for which
$\dist(X_i,\d X_{i+1})\ge\la^{i+1}ml$.
Now, we define
$\wh\Pi=\wh\Pi_{m,k}:=\cup_{i\in\N}X_i$.
Given
$x$, $x'\in\wh\Pi$,
there is
$i\in\N$
with
$x$, $x'\in X_i$.
Then, the distance
$|xx'|$
in
$\wh\Pi$
is well defined as the distance between
$x$, $x'$
in
$X_i$.
Therefore,
$\wh\Pi$
is a metric space.

\begin{pro}\label{pro:asympont} The the metric space
$\wh\Pi$
is asymptotically similar to the compact space
$\Pi$.
\end{pro}

\begin{proof} We put
$\la=\La_0=mk$
and consider a bounded subset
$A\sub\wh\Pi$
with
$\diam A\le R/\La_0$
for some
$R>1$.
There is
$i\in\N$
with
$\la^i<R\le\la^{i+1}$.
Then, any
$\la^i$-homothety
as well as
$\la^{i+1}$-homothety
is a
$\la$-quasi-homothety
with coefficient
$R$.

Because
$\dist(X_j,\d X_{j+1})\to\infty$
as
$j\to\infty$, $A$
is contained in some
$X_{j+1}$,
and we can assume that
$j\ge i$.
The problem is however that
$j$
can be much larger than
$i$,
e.g. if
$A$
is sitting near the boundary of
$X_j$.
Thus we assume that
$j\ge i$
is minimal with property
$A'\sub X_{j+1}$
for some isometric copy
$A'\sub\wh\Pi$
of
$A$.
Now, we show that
$j\le i+1$.

Assume to the contrary, that
$j\ge i+2$.
Note that
$\la=mk\ge 6$.
Then
$R/\La_0\le\la^i\le\la^{i+1}/6$.
The space
$X_{j+1}$
is the union of
$s_{m,k}^2$
subblocks each of which is isometric to
$X_{j-1}$.
Since
$j-1\ge i+1$
and the distance in
$X_{j+1}$
between any disjoint copies of
$X_{j-1}$
is at least
$\la^{i+1}$,
we obtain that
$A$
is covered by the union
$Q\sub X_{j+1}$
of copies of
$X_{j-1}$
such that there is a common point of the copies.
However, every such a union is isometric to a
subset of
$X_j$.
Hence,
$X_j$
contains an isometric copy of
$A$.
Since this contradicts our assumption on
$j$,
we conclude that
$j\le i+1$.

It follows that the
$\la^{j+1}$-homothety
$f:\Pi\to X_{j+1}\sub\wh\Pi$
is a
$\la$-quasi-homothety
with coefficient
$R$,
and its image
$X_{j+1}$
contains an isometric copy of
$A$.
Therefore,
$\wh\Pi$
is asymptotically similar to
$\Pi$.
\end{proof}

Using Theorem~\ref{thm:main2}, we obtain
$\asdim\wh\Pi=\dim\Pi=2$.
We fix an odd
$k=2l+1\ge 3$
and put
$\wh\Pi_p=\wh\Pi_{p,k}$.

\begin{cor}\label{cor:asympont} For distinct prime
$p$, $q$,
we have
$$\asdim(\wh\Pi_p\times\wh\Pi_q)<\asdim\wh\Pi_p+
  \asdim\wh\Pi_q.$$
\end{cor}

\begin{proof} If metric spaces
$X$, $X'$
are asymptotically similar to metric spaces
$Y$, $Y'$
respectively, then clearly
$X\times X'$
is asymptotically similar to
$Y\times Y'$.
Therefore,
$\asdim(\wh\Pi_p\times\wh\Pi_q)=\dim(\Pi_p\times\Pi_q)=3$,
while
$\asdim\wh\Pi_p=\dim\Pi_p=2$
and
$\asdim\wh\Pi_q=\dim\Pi_q=2$.
\end{proof}

\begin{rem}\label{rem:coarse} In \cite{Gra}, examples of
coarse spaces
$X$, $Y$
with
$\asdim(X\times Y)<\asdim X+\asdim Y$
are given. Here, the asymptotic dimension
$\asdim X$
is associated with the coarse structure
$\cE$
of
$X$
and it would more appropriate to use notation
$\cE\dim X$
for that dimension. Spaces
$X$, $Y$
are of the form
$K\times[0,1)$
where
$K$
is a classical Pontryagin surface, that is a 2-dimensional
compact space with
$\dim_{\Q}K=\dim_{\Z_q}K=1$
for every prime
$q\neq p$
and
$\dim_{\Z_p}K=2$
for some prime
$p$.
The coarse structure on
$X=K\times[0,1)$
is the topological coarse structure (in terms of the book
\cite{Ro}) or the continuously controlled coarse structure
(in terms of \cite{Gra}) induced by the compactification
$\ov X=K\times[0,1]$
of
$X$,
and
$\ov X$
is certainly metrizable.

We want to explain that these examples do not cover the
case of the classical asymptotic dimension, introduced in
\cite{Gr}, and for which we have constructed our asymptotically
self-similar Pontryagin surfaces. The usual asymptotic
dimension of a metric space
$(X,d)$
is associated with the bounded coarse structure
$\cE_d$, $\asdim X=\cE_d\dim X$,
where a set
$E\sub X\times X$
is controlled,
$E\in\cE_d$,
if and only if
$\sup\set{d(x,x')}{$(x,x')\in E$}$
is finite. However, it is known that the coarse
structure induced by any metrizable compactification
is nonmetrizable, see \cite[Example~2.53 and
Remark~2.54]{Ro}. In particular, the coarse structure
of
$X=K\times[0,1)$
above is nonmetrizable, i.e. there is no metric
$d$
on
$X$
for which
$\cE=\cE_d$,
and therefore, it does not make sense to compare
$\cE\dim X$
with the classical asymptotic dimension. Rather,
arguments from \cite[Theorem~2.5.7]{Gra} show that
$\cE\dim X=\dim K+1$
coincides with the topological dimension of
$X$, $\cE\dim X=\dim X$.
\end{rem}

\begin{tabbing}

Sergei Buyalo,\hskip11em\relax \= Nina Lebedeva,\\

{\tt sbuyalo@pdmi.ras.ru}\> {\tt lebed@pdmi.ras.ru}

\end{tabbing}

\noindent
St. Petersburg Dept. of Steklov \\
Math. Institute RAS, Fontanka 27, \\
191023 St. Petersburg, Russia

\end{document}